\documentclass[journal]{IEEEtran}
\usepackage{cite}
\usepackage{amsmath,amssymb,amsfonts}
\usepackage{algorithmic}
\usepackage{graphicx}
\usepackage{textcomp}

\usepackage{url}

\usepackage{amsmath,mathtools,amssymb,amsthm}

\usepackage[inline]{enumitem}
\usepackage{arydshln}

\DeclareMathOperator{\im}{im}

\usepackage{eso-pic}
\AddToShipoutPictureBG*{%
\AtPageUpperLeft{%
\setlength\unitlength{1in}%
\hspace*{\dimexpr0.5\paperwidth\relax}
\makebox(0,-1)[c]{
\begin{tabular}{c c}
T.~J.~{Meijer} \emph{et al.},
A Unified Non-Strict Finsler Lemma, \\
To appear in {\em IEEE Control Systems Letters}, uploaded to ArXiv \today \\
\end{tabular}}}}
\AddToShipoutPictureBG*{%
\AtPageUpperLeft{%
\setlength\unitlength{1in}%
\hspace*{\dimexpr0.5\paperwidth\relax}
\makebox(0,-21.4)[c]{
\footnotesize
\begin{tabular}{c c}
© 2024 IEEE. Personal use of this material is permitted. Permission from
IEEE must be obtained for all other uses, in any \\
current or future media, including reprinting/republishing
this material for advertising or promotional purposes, creating new \\
collective works, for resale or redistribution to servers or lists,
or reuse of any copyrighted component of this work in other works.
\end{tabular}}}}

\let\leq\relax\let\geq\relax\let\succeq\relax\let\preceq\relax
\let\hat\relax
\makeatletter
\def\hat{\@latex@warning{hat used}}
\def\leq{\@latex@warning{leq used}}\def\geq{\@latex@warning{geq used}}\def\succeq{\@latex@warning{succeq used}}\def\preceq{\@latex@warning{preceq used}}
\makeatother

\newtheorem{lemma}{Lemma}
\newtheorem{example}{Example}
\newtheorem{theorem}{Theorem}
\newtheorem{corollary}{Corollary}

\newtheorem{proposition}{Proposition}

\def\BibTeX{{\rm B\kern-.05em{\sc i\kern-.025em b}\kern-.08em T\kern-.1667em\lower.7ex\hbox{E}\kern-.125emX}}
\title{A Unified Non-Strict Finsler Lemma}

\author{T.~J.~{Meijer}, K.~J.~A.~{Scheres}, S.~{van den Eijnden}, T.~{Holicki},\\ C.~W.~{Scherer}, Fellow, IEEE, and W.~P.~M.~H.~{Heemels}, Fellow, IEEE 
\thanks{$^{1}$Tomas Meijer, Koen Scheres, Sebastiaan van den Eijnden and Maurice Heemels are with the Department of Mechanical Engineering, Eindhoven University of Technology, P.O. Box 513, 5600 MB Eindhoven, The Netherlands.
\{t.j.meijer; k.j.a.scheres; s.j.a.m.v.d.eijnden; m.heemels\}@tue.nl}%
\thanks{$^2$ Tobias Holicki and Carsten Scherer are with the Department of Mathematics, University of Stuttgart, Germany.
tobias.holicki@mathematik.uni-stuttgart.de; carsten.scherer@imng.uni-stuttgart.de}%
\thanks{This research received funding from the European Research Council (ERC) under the Advanced ERC grant agreement PROACTHIS, no. 101055384.}
\thanks{The fourth and fifth author are funded by Deutsche Forschungsgemeinschaft (DFG, German Research Foundation) under Germany’s Excellence Strategy - EXC 2075 - 390740016. They acknowledge the support by the Stuttgart Center for Simulation Science (SimTech).}
\thanks{$^{\dagger}$Corresponding author: T.~J.~{Meijer}}%
}

\begin{document}
\maketitle
\thispagestyle{empty}
\pagestyle{empty} 

\begin{abstract}
    In this paper, we present a unified general non-strict Finsler lemma. 
    This result is general in the sense that it does not impose any restrictions on the involved matrices and, thereby, it encompasses all existing non-strict versions of Finsler's lemma that do impose such restrictions.
    To further illustrate its usefulness, we showcase applications of the non-strict Finsler's lemma in deriving a structured solution to a special case of the non-strict projection lemma, and we use the unified non-strict Finsler's lemma to prove a more general version of the matrix Finsler's lemma.
\end{abstract}


\section{Introduction}
\IEEEPARstart{F}{insler's} lemma was first introduced in~\cite{Finsler1937} and is stated as follows (the adopted notation is explained at the end of this introduction).

\begin{lemma}\label{lem:sfl}
    Let $M$ and $N$ be symmetric matrices. Then, the following statements are equivalent:
    \begin{enumerate}[labelindent=5pt,labelwidth=\widthof{\ref{item:S2}},itemindent=0em,leftmargin=!,label=(S\arabic*)]
        \item\label{item:S1} there exists $\alpha\in\mathbb{R}$ such that $M+\alpha N\succ 0$;
        \item\label{item:S2} $x^\top Mx> 0$ for all $x\in\mathbb{R}^n\setminus\{0\}$ with $x^\top Nx=0$.
    \end{enumerate}
\end{lemma}
\noindent Since its introduction, Finsler's lemma (Lemma~\ref{lem:sfl}) has found application in many optimization and control problems. This includes, for example, the analysis of piecewise linear systems with sliding modes, where Lemma~\ref{lem:sfl} is used for formulating tractable conditions in the form of linear matrix inequalities (LMIs) that guarantee the decrease of a quadratic function on the sliding surface \cite{Heemels2008}. A related example comes from the stability and performance analysis of projection-based controllers, where some of the controller states are kept into constraint sets by means of projection. There, Lemma~\ref{lem:sfl} is crucial for analyzing the closed-loop system behaviour on the boundary of the constraint set, see, e.g.,~\cite{Deenen2021b,vandenEijnden2022}. In the optimization literature, Lemma~\ref{lem:sfl} is used, for instance, in the analysis of augmented Lagrangian algorithms. There, the matrix $M$ corresponds to the Hessian of the Lagrangian and $N$ typically relates to the active constraint normals. Other noteworthy applications include data-driven and model predictive control~\cite{vanWaarde2021,Nguyen2023}, robust controller synthesis~\cite{Xie1992,Oliveira2007}, and observer design~\cite{Meijer2022-nl-obs-arxiv}. All in all, Finsler's lemma has proven to be an indispensable tool in optimization and control theory. 

In its original form, Lemma~\ref{lem:sfl} deals with strict inequalities only. Hence, in the remainder of this paper, we refer to Lemma~\ref{lem:sfl} as the \emph{strict} Finsler's lemma (SFL). Considering strict inequalities only, however, may be limiting in many applications. While the implication~\ref{item:S1}$\Rightarrow$\ref{item:S2} still holds in the non-strict case, i.e., when we replace the \emph{strict} inequalities ($\succ$/$>$) by their \emph{non-strict} counterparts ($\succcurlyeq$/$\geqslant$), the converse is no longer true, as illustrated in the example below.
\begin{example}\label{ex:infeasible-case}
    Consider the matrices
        $M = \begin{bmatrix}\begin{smallmatrix}
            1 & 0\\
            0 & -1
        \end{smallmatrix}\end{bmatrix}$ and 
        $N=\begin{bmatrix}\begin{smallmatrix}
            1 & 1\\
            1 & 1
        \end{smallmatrix}\end{bmatrix}$.
    Observe that $N\succcurlyeq 0$ and, thereby, any $x\in\{\xi\in\mathbb{R}^2\mid \xi^\top N\xi=0\}$ satisfies $x\in\ker N$. For any $x\in\ker N$, we have $x=[\beta\:-\beta]^\top$ for some $\beta\in\mathbb{R}$ and thus $x^\top Mx=0$. Hence, it is clear that the non-strict version of~\ref{item:S2} holds. However, the non-strict version of~\ref{item:S1} is infeasible, as there exists no $\alpha\in\mathbb{R}$ for which $M+\alpha N\succcurlyeq 0$ since $\det M+\alpha N=-1$.
\end{example}

Although several extensions of Finsler's lemma to the non-strict case are available in the literature, their focus lies on specific settings, e.g., when $N$ is indefinite \cite[Thm.~2.3]{Moré1993} or when $N$ is positive semi-definite \cite[Thm.~4.2]{Anstreicher2000}, and their respective proofs do not allow the inclusion of the other settings straightforwardly. This may be problematic in situations where it is hard to determine \emph{a priori} whether $N$ is (in)definite, e.g., when $N$ includes decision variables in an optimization problem. Therefore, it is valuable to have a single result and proof that covers \emph{all} cases, i.e., applies to arbitrary symmetric matrices $M$, $N$. 

The main contribution of this paper is a \emph{unified} non-strict Finsler lemma (NSFL) that is \emph{general} in the sense that it encompasses all existing results and applies to any symmetric matrices $M$ and $N$, for which it provides necessary and sufficient conditions for the feasibility of the non-strict version of~\ref{item:S1}. The key challenge in doing so is to devise an additional condition that \begin{enumerate*}[label=$(\roman*)$] \item together with the non-strict version of~\ref{item:S2} implies feasibility of the non-strict version of~\ref{item:S1}, and \item is itself implied by feasibility of the non-strict version of~\ref{item:S1}.\end{enumerate*} Interestingly, the non-strict projection lemma~\cite[Thm.~1]{Meijer2024-nspl} also features such a condition, which serves as inspiration for our non-strict Finsler lemma. To establish that our main contribution truly generalizes the SFL, we show that this additional condition is trivially satisfied when~\ref{item:S2} holds. We also investigate how our unified non-strict Finsler lemma relates to several important existing results, and, in fact, we show that these existing results are special cases of our result. To further illustrate the usefulness of our results, we apply the NSFL to obtain a structured solution to a particular instance of the non-strict projection lemma~\cite{Meijer2024-nspl}. Additionally, we derive a more general version of the matrix Finsler's lemma in~\cite{vanWaarde2021} using the NSFL, and illustrate that our results are less conservative compared to those in~\cite{vanWaarde2021}. We envision that our results may be used to reduce conservatism in other contexts where non-strict inequalities are used.

The remainder of this paper is organized as follows. Section~\ref{sec:summary} gives a short survey of existing non-strict extensions of Finsler's lemma for specific cases. In Section~\ref{sec:mainresult}, we present our (unified) non-strict Finsler lemma. We showcase several applications of the non-strict Finsler's lemma in Section~\ref{sec:applications}. Finally, we provide some conclusions in Section~\ref{sec:conclusions} and all proofs can be found in the Appendix.

{\bf Notation. }The sets of real numbers, non-negative real numbers, $n$-dimensional real vectors and $n$-by-$m$ real matrices are denoted, respectively, by $
\mathbb{R}$, $\mathbb{R}_{\geqslant 0}$, $\mathbb{R}^{n}$ and $\mathbb{R}^{n\times m}$. The set of $n$-by-$n$ symmetric matrices is denoted by $\mathbb{S}^n=\{A\in\mathbb{R}^{n\times n}\mid A=A^\top\}$.
For any vectors $u\in\mathbb{R}^{p}$ and $v\in\mathbb{R}^{q}$, we denote $(u,v)=\begin{bmatrix}u^\top&v^\top\end{bmatrix}^\top$. The symbol $I$ is an identity matrix of appropriate dimension. For a symmetric matrix $S\in\mathbb{S}^n$, $S\succ 0$, $S\succcurlyeq 0$, $S\prec 0$ and $S\preccurlyeq 0$ mean that $S$ is, respectively, positive definite, positive semi-definite, negative definite,  and negative semi-definite. For a positive semi-definite matrix $S\succcurlyeq 0$, $S^\frac{1}{2}\succcurlyeq 0$ denotes its unique positive semi-definite and symmetric square root such that $S^{\frac{1}{2}}S^{\frac{1}{2}}=S$. We use the symbol $\star$ to complete a symmetric matrix, e.g.,
$\begin{bmatrix}\begin{smallmatrix} A & B\\\star & C\end{smallmatrix}\end{bmatrix}=\begin{bmatrix}\begin{smallmatrix}A & B\\B^{\top } & C\end{smallmatrix}\end{bmatrix}$. For a matrix $A\in\mathbb{R}^{n\times m}$, $\im A = \{x\in\mathbb{R}^{n}\mid x = Ay \text{ for some }y\in\mathbb{R}^{m}\}$ denotes its image, $\ker A=\{x\in\mathbb{R}^{m}\mid Ax=0\}$ its kernel,  and $A^+$ its (Moore-Penrose) pseudoinverse. The annihilator $A_{\perp}$ denotes any matrix whose columns form a basis of $\ker A$. For a subspace $\mathcal{S}\subset\mathbb{R}^{n}$, $\mathcal{S}^\perp$ denotes its complement.

\section{Existing non-strict Finsler's lemmas}\label{sec:summary}
In this section we provide a brief survey of existing non-strict extensions of Lemma~\ref{lem:sfl}. We consider two important extensions of Lemma~\ref{lem:sfl} to non-strict inequalities that consider, respectively, the case where \begin{enumerate*}[label=(\alph*)] \item $N$ is indefinite, and \item $N$ is definite\end{enumerate*}. Together, this covers all possible matrices $N$ and $M$ (because no constraints are imposed on the latter).

\subsection{The indefinite case}
First, we consider the case where $N$ is indefinite, which is covered by the well-known result found in~\cite[Thm.~2.3]{Moré1993}.
\begin{lemma}\label{lem:indef-fl}
    Let $M,N\in\mathbb{S}^{n}$ and suppose that $N$ is indefinite. Then, the following statements are equivalent:
    \begin{enumerate}[labelindent=5pt,labelwidth=\widthof{\ref{item:I2}},itemindent=0em,leftmargin=!,label=(I\arabic*),topsep=6pt]
        \item\label{item:I1} there exists $\alpha\in\mathbb{R}$ such that $M+\alpha N\succcurlyeq 0$;
        \item\label{item:I2} $x^\top Mx\geqslant 0$ for all $x\in\mathbb{R}^{n}$ with $x^\top Nx=0$.
    \end{enumerate}
\end{lemma}
\noindent The above result has proved useful in several applications such as for analyzing trust region algorithms \cite{Moré1993}. Unlike Lemma~\ref{lem:sfl}, however, Lemma~\ref{lem:indef-fl} only applies when the matrix $N$ is indefinite. According to~\cite{Moré1993}, this indefiniteness of $N$ is even necessary in the non-strict case. However, as illustrated by Lemma~\ref{lem:def-fl} in the next section, there are cases where the matrix $N$ is not indefinite but still~\ref{item:I1} and~\ref{item:I2} both hold.

\subsection{The definite case}
A non-strict extension of Lemma~\ref{lem:sfl} that does not require indefiniteness of $N$ is found in~\cite{Anstreicher2000} and repeated below.
\begin{lemma}[{\cite[Thm. 4.2]{Anstreicher2000}}]\label{lem:def-fl}
    Let $M,N\in\mathbb{S}^{n}$ be symmetric and suppose that $N\succcurlyeq 0$. Consider the following statements:
    \begin{enumerate}[labelindent=5pt,labelwidth=\widthof{\ref{item:D3}},itemindent=0em,leftmargin=!,label=(D\arabic*),topsep=6pt]
        \item\label{item:D1} there exists $\alpha\in\mathbb{R}_{\geqslant 0}$ such that $M+\alpha N\succcurlyeq 0$;
        \item\label{item:D2} $x^\top Mx\geqslant 0$ for all $x\in\mathbb{R}^{n}$ with $x^\top Nx=0$;
        \item\label{item:D3} $\ker N_\perp^\top M N_\perp=\ker N_\perp^\top M^2N_\perp~(\,\stackrel{\eqref{eq:kerBTB}}{=}\ker MN_\perp)$.
    \end{enumerate}
    If~\ref{item:D2} holds, then~\ref{item:D1} holds if and only if~\ref{item:D3} holds.
\end{lemma}
\noindent In its original publication in~\cite{Anstreicher2000}, Lemma~\ref{lem:def-fl} is stated with $N\leftarrow A^\top A$, for some matrix $A\in\mathbb{R}^{n\times m}$, and with $A_{\perp}^\top MA_\perp\succcurlyeq 0$ instead of~\ref{item:D2}. However, both formulations are equivalent. To see this, note that such $A$ exists if and only if $N\succcurlyeq 0$ and, for any $B\in\mathbb{R}^{n\times m}$, we have~\cite[Thm.~7.5.8]{Anton2003}
\begin{equation}\label{eq:kerBTB}
    \ker B^\top B = \ker B,
\end{equation}
Moreover, if $N\succcurlyeq 0$, 
\[\{\xi\in\mathbb{R}^{n}\mid x^\top Nx=0\}=\ker N\stackrel{\eqref{eq:kerBTB}}{=}\ker A.\] Hence,~\ref{item:D2} is equivalent to $A_{\perp}^\top MA_{\perp}\succcurlyeq 0$. We also note that~\cite{Anstreicher2000} only shows that \ref{item:D2} and \ref{item:D3} imply \ref{item:D1}, however, the converse also holds. To see this, observe that, since $NN_{\perp}=0$, 
\[0\preccurlyeq N_{\perp}^\top \left(M+\alpha N\right)N_{\perp} = N_{\perp}^\top MN_{\perp}.\]

While Lemma~\ref{lem:def-fl} only deals with the case where $N\succcurlyeq 0$, the case with $N\preccurlyeq 0$ follows immediately by taking $\alpha\leftarrow -\alpha$. As such, Lemma~\ref{lem:indef-fl} and~\ref{lem:def-fl} together cover all possible matrices $M$ and $N$, albeit across ``splintered'' results. To the best of our knowledge, however, no formulation in terms of non-strict inequalities is available in the literature that ``truly'' generalizes Lemma~\ref{lem:sfl} in the sense that it applies to arbitrary (symmetric) $M$ and $N$. In fact, as we will see in the next section, it turns out to be non-trivial to come up with a unified formulation and proof that deal with all cases simultaneously.

\section{Main results}
\label{sec:mainresult}
\subsection{A unified non-strict Finsler lemma}
Below, we present our unified non-strict Finsler lemma (NSFL), which forms the main result of this paper. 
\begin{theorem}\label{thm:nsfl}
  Let $M,N\in\mathbb{S}^{n}$ and consider the statements
  \begin{enumerate}[labelindent=5pt,labelwidth=\widthof{\ref{item:NS3}},itemindent=0em,leftmargin=!,label=(NS\arabic*),topsep=6pt]
    \item\label{item:NS1} there exists $\alpha\in\mathbb{R}$ such that $M+\alpha N\succcurlyeq 0$;
    \item\label{item:NS2} $x^\top Mx\geqslant 0$ for all $x\in\mathbb{R}^{n}$ with $x^\top Nx=0$;
    \item\label{item:NS3} $\ker N\cap \{\xi\in\mathbb{R}^n\mid \xi^\top M\xi=0\}\subset \ker M$.
  \end{enumerate}
  Then,~\ref{item:NS1} holds if and only if~\ref{item:NS2} and~\ref{item:NS3} hold. If, in addition to~\ref{item:NS2} and~\ref{item:NS3}, it holds that $N\succcurlyeq 0$ (or $N\preccurlyeq 0$), then~\ref{item:NS1} holds for some $\alpha>0$ (or $\alpha <0$, respectively).
\end{theorem}
\noindent Since Theorem~\ref{thm:nsfl} does not involve any additional assumption on the matrix $N$, it provides genuine necessary and sufficient conditions for feasibility of the LMI problem in \ref{item:NS1}. In contrast to \ref{item:NS1} and~\ref{item:NS2}, which can be equivalently formulated as the LMI $N_{\perp}^\top MN_{\perp}\succcurlyeq 0$, the condition in \ref{item:NS3} does not correspond to an LMI. Instead, it is a nontrivial coupling condition on the matrices $M$ and $N$ that is absent in both the SFL (Lemma~\ref{lem:sfl}) and Lemma~\ref{lem:indef-fl}. Interestingly, Lemma~\ref{lem:def-fl} does feature a type of coupling condition in~\ref{item:D3}, whose relation to~\ref{item:NS3} we investigate in more detail in Section~\ref{sec:discussion}.
The condition in \ref{item:NS3} is reminiscent of the coupling condition that is found in the non-strict projection lemma~\cite{Meijer2024-nspl}. Often,~\ref{item:NS3} is satisfied inherently as a consequence of the underlying problem structure. This is illustrated in several applications in the remainder of this and the next section and the respective proofs. We conclude this section by revisiting Example~\ref{ex:infeasible-case}, for which we found that~\ref{item:NS1} was infeasible and, thus, by Theorem~\ref{thm:nsfl},~\ref{item:NS3} must be false, as confirmed below.
\begin{example}
    As~\ref{item:NS1} is not feasible for $M$ and $N$ in Example~\ref{ex:infeasible-case} even though~\ref{item:NS2} holds, by Theorem~\ref{thm:nsfl},~\ref{item:NS3} cannot hold. Indeed, we find 
    \[\ker N\cap \{\xi\in\mathbb{R}^{2}\mid \xi^\top M\xi=0\}=\im{} (1,-1),\]
    which is not contained in $\ker M=\{0\}$. 
\end{example}

\subsection{Discussion}\label{sec:discussion}
Next, we elaborate on the relation between Theorem~\ref{thm:nsfl} and the existing results, i.e., Lemmas~\ref{lem:sfl}-\ref{lem:def-fl}. For each of these lemmas, the most difficult step in the proof is showing that the LMI-based condition holds, i.e., that~\ref{item:S2} implies~\ref{item:S1},~\ref{item:I2} implies~\ref{item:I1}, and~\ref{item:D2} (with~\ref{item:D3}) implies~\ref{item:D1}. In this section, we illustrate that these implications are all consequences of Theorem~\ref{thm:nsfl}. In doing so, we also provide several examples illustrating how the non-trivial coupling condition~\ref{item:NS3} can be verified.

    First, we show that the coupling condition \ref{item:NS3} is trivially satisfied if \ref{item:NS2} is strict, i.e., if \ref{item:S2} is satisfied.
    \begin{proposition}\label{prop:special-case-sfl}
        Let $M,N\in\mathbb{S}^n$ and suppose that~\ref{item:S2} holds. Then,~\ref{item:NS3} also holds.
    \end{proposition}
    \noindent By a similar perturbation argument as in~\cite[Cor.~1]{Meijer2024-nspl}, we can apply Theorem~\ref{thm:nsfl} to show the equivalence in Lemma~\ref{lem:sfl}.

Next, we consider the non-strict case, where, interestingly, indefiniteness of $N$ together with~\ref{item:I2} implies~\ref{item:NS3}.
\begin{proposition}\label{prop:special-case-indef}
    Let $M,N\in\mathbb{S}^n$ and suppose that $N$ is indefinite and that~\ref{item:I2} holds. Then,~\ref{item:NS3} also holds.
\end{proposition}
\noindent Using Theorem~\ref{thm:nsfl}, we can immediately conclude also that, if~\ref{item:I2} holds, then~\ref{item:I1} must hold. In other words, Theorem~\ref{thm:nsfl} is at least as general as Lemma~\ref{lem:indef-fl}.
It also turns out that~\ref{item:D3} is equivalent to \ref{item:NS3}, if \ref{item:D2} holds, regardless of positive/negative semi-definiteness of $N$.
\begin{proposition}\label{prop:equivalence-coupling}
    Let $M,N\in\mathbb{S}^n$ and suppose that~\ref{item:D2} holds. Then, the statements~\ref{item:D3} and~\ref{item:NS3} are equivalent.
\end{proposition}
\noindent Proposition~\ref{prop:equivalence-coupling} shows that the coupling condition~\ref{item:D3} turns out to be stronger than initially shown in~\cite{Anstreicher2000} and that Theorem~\ref{thm:nsfl} is at least as general as Lemma~\ref{lem:def-fl}. Proposition~\ref{prop:equivalence-coupling} also means that \ref{item:D3} is implied by~\ref{item:D1} and, thus,~\ref{item:D3} can be used to verify \ref{item:NS3}, which leads to the corollary below.
\begin{corollary}\label{cor:ns3-replacement}
    Let $M,N\in\mathbb{S}^n$. Then, \ref{item:NS1} holds if and only if \ref{item:NS2} and~\ref{item:D3} hold.
\end{corollary}

\section{Applications}\label{sec:applications}
\subsection{Structured solution to the non-strict projection lemma}
Interestingly, Theorem~\ref{thm:nsfl} can be used to obtain a structured solution to the non-strict projection lemma (NSPL), which was recently presented in~\cite{Meijer2024-nspl}. While this is generally not possible, we provide a structured solution below for the case where either $V=I$ or $U=I$ (in the notation of~\cite{Meijer2024-nspl}).
\begin{lemma}\label{lem:structured-nspl}
    Let $U\in\mathbb{R}^{m\times n}$ and $Q\in\mathbb{S}^{n}$. There exists $X\in\mathbb{R}^{m\times n}$ satisfying
    \begin{equation}
        Q+U^\top X + X^\top U\succcurlyeq 0,\label{eq:QUX}
    \end{equation}
    if and only if
    \begin{gather}
        U_{\perp}^\top QU_{\perp}\succcurlyeq 0,\label{eq:UQU}        
    \end{gather}
    If, in addition, it holds that \begin{equation}\ker U\cap \{\xi\in\mathbb{R}^{n}\mid \xi^\top Q\xi=0\}\subset\ker Q,\label{eq:one-sided-coupling}\end{equation}
    then there exists $\alpha>0$ such that~\eqref{eq:QUX} holds with
    \begin{equation}
        X = \alpha U.\label{eq:structured-X}
    \end{equation}
\end{lemma}
\noindent Note that the equivalence of~\eqref{eq:QUX} and~\eqref{eq:UQU} follows immediately from~\cite[Thm.~1]{Meijer2024-nspl}, however, the added value in Lemma~\ref{lem:structured-nspl} is in the provided structured solution~\eqref{eq:structured-X} for the case where~\eqref{eq:one-sided-coupling} holds. The crucial step in obtaining this structured solution consists of applying Theorem~\ref{thm:nsfl} with $N=U^\top U$. The structured solution in~\eqref{eq:structured-X} can, depending on the structure of $U$ itself, often be used to simplify the matrix inequality~\eqref{eq:QUX}.

\subsection{A matrix Finsler lemma}
In this section, we derive, based on Theorem~\ref{thm:nsfl}, a generalization of the matrix Finsler's lemma~\cite[Thm.~1]{vanWaarde2021}. 
\begin{lemma}\label{lem:MFL}
    Let 
    \begin{equation*}
        M = \begin{bmatrix}
            M_{11} & M_{12}\\
            M_{12}^\top & M_{22}
        \end{bmatrix}\in\mathbb{S}^{n+m}\text{ and }N = \begin{bmatrix}
            N_{11} & N_{12}\\
            N_{12}^\top & N_{22}
        \end{bmatrix}\in\mathbb{S}^{n+m}
    \end{equation*}
    with $N_{22}\preccurlyeq 0$, $N_{11}-N_{12}N_{22}^+N_{12}^\top = 0$, $M_{22}\preccurlyeq 0$ 
 and $(I-N_{22}N_{22}^+)N_{12}^\top=0$. Consider
    \begin{enumerate}[labelindent=5pt,labelwidth=\widthof{\ref{item:M2}},itemindent=0em,leftmargin=!,label=(M\arabic*)]\addtocounter{enumi}{1}
        \item\label{item:M2} $\begin{bmatrix}
            I\\
            Z
        \end{bmatrix}^\top M\begin{bmatrix}
            I\\
            Z
        \end{bmatrix}\succcurlyeq 0$ for all $Z$ with $\begin{bmatrix}
            I\\
            Z
        \end{bmatrix}^\top N\begin{bmatrix}
            I\\
            Z
        \end{bmatrix}=0$.
    \end{enumerate}
    Then,~\ref{item:NS1} holds if and only if~\ref{item:M2} and~\ref{item:NS3} hold.
\end{lemma}
\noindent Observe that the matrix inequality 
\begin{equation*}
    \begin{bmatrix}
        I\\
        Z
    \end{bmatrix}^\top M\begin{bmatrix}
        I\\
        Z
    \end{bmatrix}\succcurlyeq 0,
\end{equation*}
is non-strict. Similar results have been published in~\cite[Thm.~1]{vanWaarde2021} and~\cite[Thm.~4.8]{vanWaarde2023b}. In contrast with~\cite[Thm.~1]{vanWaarde2021}, we do not require $M_{12}=0$, but allow a larger set of matrices $M$ to be considered. Note also that $(I-N_{22}N_{22}^+)N_{12}^\top =0$ is equivalent to the existence of $G$ such that $N_{22}G=N_{12}^\top$. The latter is assumed in~\cite[Thm.~1]{vanWaarde2021} along with stronger assumptions on $G$ and $M$. This is not required here and all assumptions are directly formulated in terms of the matrices $M$ and $N$. Instead, we exploit~\ref{item:NS3} to prove similar properties, however,~\ref{item:NS3} is part of the necessary and sufficient conditions and, thereby, does not impose additional conservatism. Moreover, these observations also mean that~\cite[Thm.~1]{vanWaarde2021} can be recovered as a special case of Lemma~\ref{lem:MFL}. To emphasize that these observations render Lemma~\ref{lem:MFL} more general than~\cite[Thm.~4.8]{vanWaarde2023b}, we consider the numerical example below.

\begin{example}\label{ex:0110}
    Consider the matrices
    \begin{equation*}\def\arraystretch{1.2}
        M = \left[\begin{array}{@{}c;{2pt/2pt}c@{}}
            M_{11} & M_{12}\\\hdashline[2pt/2pt]
            M_{12}^\top & M_{22}
        \end{array}\right] = \left[\begin{array}{@{}c;{2pt/2pt}c@{}}
            2 & 1\\\hdashline[2pt/2pt]
            1 & 0
        \end{array}\right]
    \end{equation*}
    and
    \begin{equation*}\def\arraystretch{1.2}
        N = \left[\begin{array}{@{}c;{2pt/2pt}c@{}}
            N_{11} & N_{12}\\\hdashline[2pt/2pt]
            N_{12}^\top & N_{22}
        \end{array}\right] = \left[\begin{array}{@{}c;{2pt/2pt}c@{}}
            0 & 0\\\hdashline[2pt/2pt]
            0 & -1
        \end{array}\right].
    \end{equation*}
    Clearly, $\ker M_{22}=\mathbb{R}\not\subset \ker M_{12}=\{0\}$ so that~\cite[Thm.~1]{vanWaarde2021} and \cite[Thm.~4.8]{vanWaarde2023b} cannot be applied here. However, we can still use Lemma~\ref{lem:MFL} to conclude that there exists $\alpha\in\mathbb{R}$ such that $M+\alpha N\succcurlyeq 0$. To see this,
    note that $N_{22}\preccurlyeq0$, $N_{11}-N_{12}N_{22}^+N_{12}^\top=0$, $M_{22}\preccurlyeq0$ and $(I-N_{22}N_{22}^+)N_{12}^\top=0$. Moreover,
    \begin{equation*}
        \begin{bmatrix}
            1\\
            z
        \end{bmatrix}^\top N\begin{bmatrix}
            1\\
            z
        \end{bmatrix}=0
    \end{equation*}
    implies that $z=0$ and, thus, 
    \begin{equation*}
        \begin{bmatrix}
            1\\
            z
        \end{bmatrix}^\top M\begin{bmatrix}
            1\\
            z
        \end{bmatrix}=2\geqslant 0.
    \end{equation*}
    Hence,~\ref{item:M2} holds. Finally, we have $\ker N = \im{} (1,0)$ and $\ker N\cap \{\xi\in\mathbb{R}^2\mid\xi^\top M\xi=0\}=\{0\}\subset\ker M$, thus,~\ref{item:NS3} holds. Using Lemma~\ref{lem:MFL}, we conclude that, indeed,~\ref{item:NS1} holds.
\end{example}

The unified matrix S-lemma and Finsler's lemma in~\cite[Cor.~4.13]{vanWaarde2023b}, on the other hand, requires that $M_{11}+\alpha N_{11}\succ 0$ for some $\alpha\in\mathbb{R}$, due to the the parameter $\beta>0$ in~\cite[Cor.~4.13]{vanWaarde2023b}. The example below illustrates that this may also be conservative.
\begin{example}
    Consider the matrices
    \begin{equation*}\def\arraystretch{1.2}
        M = \left[\begin{array}{@{}c;{2pt/2pt}c@{}}
            M_{11} & M_{12}\\\hdashline[2pt/2pt]
            M_{12}^\top & M_{22}
        \end{array}\right] = \left[\begin{array}{@{}c;{2pt/2pt}c@{}}
            0 & 0\\\hdashline[2pt/2pt]
            0 & -1
        \end{array}\right]
    \end{equation*}
    and
    \begin{equation*}\def\arraystretch{1.2}
        N = \left[\begin{array}{@{}c;{2pt/2pt}c@{}}
            N_{11} & N_{12}\\\hdashline[2pt/2pt]
            N_{12}^\top & N_{22}
        \end{array}\right] = \left[\begin{array}{@{}c;{2pt/2pt}c@{}}
            0 & 0\\\hdashline[2pt/2pt]
            0 & -1
        \end{array}\right].
    \end{equation*}
    Clearly, $M_{11}+\alpha N_{11}$ is never positive definite, so~\cite[Cor.~4.13]{vanWaarde2023b} does not apply. However, we can still use Lemma~\ref{lem:MFL} to conclude that there exists $\alpha\in\mathbb{R}$ such that $M+\alpha N\succcurlyeq 0$. To see this, 
    note that $N_{22}\preccurlyeq0$, $N_{11}-N_{12}N_{22}^+N_{12}^\top=0$, $M_{22}\preccurlyeq0$ and $(I-N_{22}N_{22}^+)N_{12}^\top=0$. Moreover,
    \begin{equation*}
        \begin{bmatrix}
            1\\
            z
        \end{bmatrix}^\top N\begin{bmatrix}
            1\\
            z
        \end{bmatrix}=0
    \end{equation*}
    implies that $z=0$ and, thus, 
    \begin{equation*}
        \begin{bmatrix}
            1\\
            z
        \end{bmatrix}^\top M\begin{bmatrix}
            1\\
            z
        \end{bmatrix}=0\geqslant 0.
    \end{equation*}
    Hence,~\ref{item:M2} holds. Lastly, we have $\ker N = \im{} (1,0)$ and $\ker N\cap \{\xi\in\mathbb{R}^2\mid\xi^\top M\xi=0\}=\im{} (1,0)=\ker M$, thus,~\ref{item:NS3} holds. Using Lemma~\ref{lem:MFL}, we conclude that, indeed,~\ref{item:NS1} holds.
\end{example}

The above examples illustrate how, by using the necessary and sufficient version of the non-strict Finsler's lemma, derivatives such as the matrix Finsler's lemma can be obtained that do not impose (additional) restrictive conditions on the matrices $M$ and $N$.

\section{Conclusions}\label{sec:conclusions}
In this paper, we presented a unified non-strict Finsler lemma that encompasses important existing non-strict extensions of Finsler's lemma. It turns out that, in order to do so, an additional non-trivial coupling condition is needed that, combined with the standard (but now non-strict) conditions of Finsler's lemma, provides necessary and sufficient conditions for arbitrary (symmetric) matrices. We showed that this coupling condition is trivially satisfied in the context of the original strict Finsler's lemma and that our novel non-strict Finsler lemma generalizes the existing results. Finally, we used this non-strict Finsler lemma to construct a structured solution for a specific instance of the non-strict projection lemma, and to prove a more general matrix Finsler lemma, for which we show by example that it is less conservative than existing versions. We envision that the necessary and sufficient conditions provided in this paper provide insights in and guidelines for reducing conservatism in other non-strict inequalities that are useful, e.g., in the context of data-driven control.




\bibliographystyle{IEEEtran}        
\bibliography{phd-bibtex-local-copy}           

\begin{thebibliography}{10}
\providecommand{\url}[1]{#1}
\csname url@samestyle\endcsname
\providecommand{\newblock}{\relax}
\providecommand{\bibinfo}[2]{#2}
\providecommand{\BIBentrySTDinterwordspacing}{\spaceskip=0pt\relax}
\providecommand{\BIBentryALTinterwordstretchfactor}{4}
\providecommand{\BIBentryALTinterwordspacing}{\spaceskip=\fontdimen2\font plus
\BIBentryALTinterwordstretchfactor\fontdimen3\font minus \fontdimen4\font\relax}
\providecommand{\BIBforeignlanguage}[2]{{%
\expandafter\ifx\csname l@#1\endcsname\relax
\typeout{** WARNING: IEEEtran.bst: No hyphenation pattern has been}%
\typeout{** loaded for the language `#1'. Using the pattern for}%
\typeout{** the default language instead.}%
\else
\language=\csname l@#1\endcsname
\fi
#2}}
\providecommand{\BIBdecl}{\relax}
\BIBdecl

\bibitem{Finsler1937}
P.~{Finsler}, ``Über das {V}orkommen definiter und semidefiniter {F}ormen in {S}charen quadratischer {F}ormen,'' \emph{Comment. Math. Helv.}, vol.~9, pp. 188--192, 1936--1937.

\bibitem{Heemels2008}
W.~P.~M.~H. {Heemels} and S.~{Weiland}, ``Input-to-state stability and interconnection of discontinuous dynamical systems,'' \emph{Automatica}, vol.~44, pp. 3079--3086, 2008.

\bibitem{Deenen2021b}
D.~A. {Deenen}, B.~{Sharif}, S.~J.~A.~M. {van den Eijnden}, H.~{Nijmeijer}, W.~P.~M.~H. {Heemels}, and M.~{Heertjes}, ``Projection-based integrators for improved motion control: {Formalization}, well-posedness and stability of hybrid integrator-gain systems,'' \emph{Automatica}, vol. 113, p. 109830, 2021.

\bibitem{vandenEijnden2022}
S.~J.~A.~M. {van den Eijnden}, W.~P.~M.~H. {Heemels}, H.~{Nijmeijer}, and M.~{Heertjes}, ``Stability and performance analysis of hybrid integrator–gain systems: {A} linear matrix inequality approach,'' \emph{Nonlinear Anal.: Hybrid Syst.}, vol.~45, p. 101192, 2022.

\bibitem{vanWaarde2021}
H.~J. {van Waarde} and M.~K. {Camlibel}, ``A matrix {Finsler}'s lemma with applications to data-driven control,'' in \emph{60th IEEE Conf. Decis. Control}, 2021, pp. 5777--5782.

\bibitem{Nguyen2023}
H.~H. Nguyen, M.~Friedel, and R.~Findeisen, ``{LMI}-based data-driven robust model predictive control,'' \emph{IFAC-PapersOnLine}, vol.~56, no.~2, pp. 4783--4788, 2023, 22nd IFAC World Congress.

\bibitem{Xie1992}
L.~Xie, M.~Fu, and C.~de~Souza, ``{$H_\infty$} control and quadratic stabilization of systems with parameter uncertainty via output feedback,'' \emph{IEEE Trans. Autom. Control}, vol.~37, no.~8, pp. 1253--1256, 1992.

\bibitem{Oliveira2007}
R.~C. Oliveira and P.~L. Peres, ``Parameter-dependent {LMI}s in robust analysis: Characterization of homogeneous polynomially parameter-dependent solutions via {LMI} relaxations,'' \emph{IEEE Trans. Autom. Control}, vol.~52, no.~7, pp. 1334--1340, 2007.

\bibitem{Meijer2022-nl-obs-arxiv}
T.~J. {Meijer}, V.~S. {Dolk}, M.~S. {Chong}, and W.~P.~M.~H. {Heemels}, ``Robust observer design for polytopic discrete-time nonlinear descriptor systems,'' 2022, preprint: \url{https://arxiv.org/abs/2207.04290}.

\bibitem{Moré1993}
J.~J. {Mor\'{e}}, ``Generalizations of the trust region problem,'' \emph{Optim. Methods Softw.}, vol.~2, pp. 189--209, 1993.

\bibitem{Anstreicher2000}
K.~M. Anstreicher and M.~H. Wright, ``A note on the augmented {H}essian when the reduced {H}essian is semidefinite,'' \emph{SIAM Journal on Optimization}, vol.~11, no.~1, pp. 243--253, 2000.

\bibitem{Meijer2024-nspl}
T.~J. {Meijer}, T.~{Holicki}, S.~J.~A.~M. {van den Eijnden}, C.~W. {Scherer}, and W.~P.~M.~H. {Heemels}, ``The non-strict projection lemma,'' \emph{IEEE Trans. Autom. Control}, pp. 1--8, 2024.

\bibitem{Anton2003}
H.~{Anton} and R.~C. {Busby}, \emph{Contemporary Linear Algebra}.\hskip 1em plus 0.5em minus 0.4em\relax John Wiley \& Sons, Inc., 2003.

\bibitem{vanWaarde2023b}
H.~J. {van Waarde}, M.~K. {Camlibel}, J.~{Eising}, and H.~L. {Trentelman}, ``Quadratic matrix inequalities with applications to data-based control,'' \emph{SIAM J. Control Optim.}, vol.~61, no.~4, pp. 2251--2281, 2023.

\bibitem{Boyd1994}
S.~{Boyd}, L.~{El Ghaoui}, E.~{Feron}, and V.~{Balakrishnan}, \emph{Linear matrix inequalities in control theory}, ser. Studies in Applied Mathematics.\hskip 1em plus 0.5em minus 0.4em\relax SIAM, 1994, vol.~15.

\end{thebibliography}
                                
\section*{Appendix}
\subsection{Proof of Theorem~\ref{thm:nsfl}}
{\bf Necessity: } Suppose that~\ref{item:NS1} holds. Then,~\ref{item:NS2} follows by substituting $x^\top Nx=0$. Next, we show that~\ref{item:NS3} holds as well. Let $x\in\ker N\cap\{\xi\in\mathbb{R}^n\mid \xi^\top M\xi=0\}$ and $S\in\mathbb{S}^{n}$ such that $S^2=M+\alpha N$, then
    \[\|Sx\|^2 = x^\top (M+\alpha N)x = x^\top Mx = 0.\]
It follows that
    \[0 = Sx = S^2x = (M+\alpha N)x = Mx,\]
whereby~\ref{item:NS3} holds.

{\bf Sufficiency: } 
Suppose that~\ref{item:NS2} and~\ref{item:NS3} hold. Let $T\in\mathbb{R}^{n\times n}$ be a non-singular matrix, whose columns in the partition $T=\begin{bmatrix} T_1 & T_2 & T_3\end{bmatrix}$ are chosen to satisfy
\begin{align}
    \im T_2 &= \ker N\cap \ker M,\label{eq:T2}\\
    \im \begin{bmatrix}
        T_1 & T_2
    \end{bmatrix} &= \ker N.\label{eq:T12}
\end{align}
It follows that~\ref{item:NS1} is equivalent to
\begin{equation}\label{eq:Y}
    Y\coloneqq T^\top (M+\alpha N)T \succcurlyeq 0.
\end{equation}
    Using~\eqref{eq:T12} and~\eqref{eq:T2}, we partition $W\coloneqq T^\top MT$ and $V\coloneqq T^\top NT$ in accordance with $T$ to obtain
\begin{equation*}\def\arraystretch{1.2}
    W = \left[\begin{array}{@{}c;{2pt/2pt}c;{2pt/2pt}c@{}}
        W_{11} & 0 & W_{13}\\\hdashline[2pt/2pt]
        0 & 0 & 0\\\hdashline[2pt/2pt]
        W_{13}^\top & 0 & W_{33}
    \end{array}\right]\text{ and }V = \left[\begin{array}{@{}c;{2pt/2pt}c;{2pt/2pt}c@{}}
        0 & 0 & 0\\\hdashline[2pt/2pt]
        0 & 0 & 0\\\hdashline[2pt/2pt]
        0 & 0 & V_{33}
    \end{array}\right].
\end{equation*}
Due to~\ref{item:NS2}, we have $W_{11} = T_1^\top MT_1 \succcurlyeq 0$, and, by~\ref{item:NS3}, 
\begin{equation}\label{eq:W11pd}
    W_{11} \succ 0.
\end{equation}
To see this, take $x=(x_1,0,0)^\top$ with $x_1$ such that $x_1^\top W_{11}x_1=0$. Observe that $x^\top Wx=0$ and $Vx=0$, i.e.,
\[x^\top\in\ker V\cap \{\xi\in\mathbb{R}^{n}\mid \xi^\top M\xi=0\}.\]
By~\ref{item:NS3}, we infer that $x\in\ker W$. It follows that $x\in\ker W\cap\ker V$ and, hence, $x_1=0$, whereby $W_{11}\succ 0$.

Using the non-strict Schur complement~\cite[pp.~8, 28]{Boyd1994}, $Y\succcurlyeq 0$ if and only if
\begin{equation}
    W_{33} - W_{13}^\top W_{11}^{-1}W_{13} + \alpha V_{33} \succcurlyeq 0.\label{eq:schured}
\end{equation}
We now prove that there exists $\alpha\in\mathbb{R}$ satisfying \eqref{eq:schured}. To this end, let $Q\coloneqq W_{33}-W_{13}^\top W_{11}^{-1}W_{13}$ and $Q_n\coloneqq Q + (1/n)I$ for $n\in\mathbb{N}$. Observe that, using~\eqref{eq:T12}, $V_{33}$ is full rank, and that, by~\eqref{eq:schured}, $Q_n$ satisfies~\ref{item:S2} for all $n\in\mathbb{N}$. Hence, by Lemma~\ref{lem:sfl}, for any $n\in\mathbb{N}$, there exists $\alpha(n)\in\mathbb{R}$ such that $Q_n+\alpha(n)V_{33}\succ 0$. We now distinguish three cases: Firstly, if $\alpha(n)\rightarrow \infty$ as $n\rightarrow\infty$, then it follows from $\alpha(n)V_{33}\succ -Q_n$ that $V_{33}\succcurlyeq 0$. Since $V_{33}$ is also full rank, we have $V_{33}\succ 0$. Hence, any $\alpha\geqslant -\lambda_{\min}(Q)/\lambda_{\min}(V_{33})$ satisfies~\eqref{eq:schured}. Secondly and similarly, if $\alpha(n)\rightarrow -\infty$ as $n\rightarrow\infty$, then $V_{33}\prec 0$. Hence, any $\alpha \leqslant -\lambda_{\min}(Q)/\lambda_{\min}(V_{33})$ satisfies~\eqref{eq:schured}. Thirdly, if neither $\alpha(n)\rightarrow \infty$ nor $\alpha(n)\rightarrow-\infty$ as $n\rightarrow\infty$, then $\{\alpha(n)\}_{n\rightarrow \infty}$ is bounded so a subsequence converges to some $\alpha\in\mathbb{R}$, and, as $Q_n+\alpha(n)V_{33} \succ 0$, $Q+\alpha V_{33}\succcurlyeq 0$. Hence, \ref{item:NS1} holds. Finally, we note that, if $N\succcurlyeq 0$, then~\ref{item:NS1} holds with $\alpha\leftarrow \bar{\alpha}$ for any $\bar{\alpha}\geqslant \min\{0,\alpha\}$. Similarly, if $N\preccurlyeq 0$, then we can take any $\bar{\alpha}\leqslant \max\{0,\alpha\}$. In the above, we assumed that $T_1$ is nonempty. If $T_1$ is empty, $W_{11}$ and $W_{13}$ are empty and \eqref{eq:Y} is equivalent to $W_{33} + \alpha V_{33} \succcurlyeq 0$, for which the proof holds \emph{mutatis mutandis}. \hfill\qed

\subsection{Proof of Proposition~\ref{prop:special-case-sfl}}
By~\ref{item:S2}, $x^\top Mx>0$ for all $x\in\ker N\setminus\{0\}$. It follows that $\ker N\cap\{\xi\in\mathbb{R}^n\mid \xi^\top M\xi=0\}=\{0\}\subset\ker M$. \hfill\qed

\subsection{Proof of Proposition~\ref{prop:special-case-indef}}
We proceed by contradiction. To this end, suppose that~\ref{item:I2} holds (and $N$ is indefinite), but that~\ref{item:NS3} does \emph{not} hold. Then, there exists $x\in\ker N$ with $x^\top Mx=0$ but $z\coloneqq Mx\neq 0$. We can find some $y$ with $y^\top Ny = 0$ and $x^\top  My\neq 0$. If $z$ satisfies $z^\top Nz = 0$, we can simply take $y=z$ since $y^\top Mx = \|z\|^2 >0$. If $z$ satisfies $z^\top Nz >0$, however, we can find, by indefiniteness of $N$, some $d$ with $d^\top Nd <0$ and $d^\top z\geqslant 0$ (note that, if $d^\top z<0$, we can take $-d$). Consider the quadratic function 
\[ f(t)\coloneqq (z+t(d-z))^\top N(z+t(d-z)).\]
Since $f(0)=z^\top Nz>0$ and $f(1) = d^\top Nd<0$, $f$ has a root $t_0\in(0,1)$. Since $d^\top z\geqslant 0$ and $t_0<1$, we infer 
\[(z+t_0(d-z))^\top z = (1-t_0)\|z\|^2 + t_0d^\top z>0\]
and take $y=z+t_0(d-z)$, which satisfies $y^\top Mx=y^\top z>0$ and $y^\top Ny=f(t_0)=0$. Similar arguments apply if $z^\top Nz<0$.

Using $y^\top Ny=0$ and $x\in\ker N$, we find 
\[(ty + x)^\top N(ty + x) = t^2y^\top Ny + 2ty^\top Nx + x^\top Nx = 0\]
for all $t\in\mathbb{R}$. Then, due to~\ref{item:I2}, we find that
\[g(t) \coloneqq (ty + x)^\top M(ty + x) = t^2y^\top My + 2ty^\top Mx \geqslant 0\]
for all $t\in\mathbb{R}$. However, since $g(0)=0$ and 
\[g^\prime (0) = 2y^\top Mx\neq 0,\]
$g(t)<0$ for some $t\in\mathbb{R}$, which is a contradiction. Thus,~\ref{item:NS3} holds. \hfill\qed

\subsection{Proof of Proposition~\ref{prop:equivalence-coupling}}
    {\bf \ref{item:D3}$\Rightarrow$\ref{item:NS3}:} Suppose that~\ref{item:D3} holds. Let $x\in\ker N$ be such that $x^\top Mx=0$. It follows that $x=N_{\perp}y$ for some $y$, so that $y^\top N_{\perp}^\top MN_{\perp}y=0$. By~\ref{item:D2}, $N_{\perp}^\top MN_{\perp}\succcurlyeq 0$. It follows that $\|(N_{\perp}^\top MN_{\perp})^{\frac{1}{2}}y\|^2 =0$ and, thereby, $N_{\perp}^\top MN_{\perp}y=0$. Thus, $y\in\ker N_\perp^\top MN_\perp$, which, by \ref{item:D3}, implies that $y\in\ker N_\perp^\top M^2N_\perp$. Using~\eqref{eq:kerBTB}, we conclude that $y\in\ker MN_\perp$. Since $x=N_\perp y$, we have $x\in\ker M$, whereby \ref{item:NS3} holds.
    
    {\bf \ref{item:D3}$\Leftarrow$\ref{item:NS3}:} Suppose that~\ref{item:NS3} holds. First, we show $\ker N_{\perp}^\top M^2N_{\perp} \subset \ker N_{\perp}^\top MN_{\perp}.$ Let $x\in\ker N_{\perp}^\top M^2N_{\perp}$. It follows, using $M=M^\top$ and~\eqref{eq:kerBTB}, that $0=\|N_{\perp}^\top M^2N_{\perp}x\| = \|MN_{\perp}x\|$ and, thus, $N_{\perp}^\top MN_\perp x=0$. We conclude that $\ker N_{\perp}^\top M^2N_{\perp} \subset\ker N_{\perp}^\top MN_{\perp}$. Next, we show $\ker N_{\perp}^\top M^2N_{\perp} \supset \ker N_{\perp}^\top MN_{\perp}.$ Let $x\in\ker N_{\perp}^\top MN_{\perp}$. Then, $0=x^\top N_{\perp}^\top MN_{\perp}x = y^\top My$ with $y=N_{\perp}x$. Since $y\in\ker N$ and $y^\top My=0$, it holds, using~\ref{item:NS3}, that $y\in\ker M$ and, thus, $x\in\ker MN_{\perp}$. Finally, using $M=M^\top$ and~\eqref{eq:kerBTB}, we conclude that $x\in\ker N_{\perp}^\top M^2N_{\perp}$, whereby $\ker N_{\perp}^\top M^2N_{\perp} \supset \ker N_{\perp}^\top MN_{\perp}$ holds.\hfill\qed

\subsection{Proof of Lemma~\ref{lem:structured-nspl}}
Since the equivalence in Lemma~\ref{lem:structured-nspl} follows immediately from the non-strict projection lemma~\cite[Thm.~1]{Meijer2024-nspl} with $V=I$, we focus on constructing the structured solution in~\eqref{eq:structured-X}. Suppose that~\eqref{eq:UQU} and~\eqref{eq:one-sided-coupling} hold and let $M=Q$ and $N=U^\top U$. It follows, using that $\ker N=\ker U$ due to~\eqref{eq:kerBTB}, that
\begin{equation*}
    \ker N\cap \{\xi\in\mathbb{R}^n\mid\xi^\top M\xi=0\}\stackrel{\eqref{eq:one-sided-coupling}}{\subset}\ker Q = \ker M,
\end{equation*}
whereby~\ref{item:NS3} holds. Next, we show that~\ref{item:NS2} holds. To this end, let $x\in\mathbb{R}^{n}$ be such that $x^\top Nx=0$. It follows that $x^\top U^\top Ux = \|Ux\|^2 = 0$ and, thus, $Ux=0$, i.e., $x\in\ker N$. Then, there exists $\eta$ such that $x=U_{\perp}\eta$ and we find that
\begin{equation*}
    x^\top Mx = \eta^\top U_{\perp}^\top QU_{\perp} \eta\stackrel{\eqref{eq:UQU}}{\geqslant} 0,
\end{equation*}
which shows that, indeed,~\ref{item:NS2} holds. Hence, by application of Theorem~\ref{thm:nsfl} and $N\succcurlyeq 0$, there exists $\alpha>0$ such that
\begin{equation*}
    M+2\alpha N = Q+2\alpha U^\top U\succcurlyeq 0,
\end{equation*}
which is~\eqref{eq:QUX} with $X = \alpha U$ and completes our proof.\hfill\qed


\subsection{Proof of Lemma~\ref{lem:MFL}}
{\bf Necessity:} Suppose that~\ref{item:NS1} holds. It trivially follows that, for all $Z\in\mathbb{R}^{m\times n}$ with $\begin{bmatrix} I & Z^\top\end{bmatrix} N\begin{bmatrix} I & Z^\top\end{bmatrix}^\top =0$, 
\[
    \begin{bmatrix} I \\ Z\end{bmatrix}^\top M\begin{bmatrix} I \\ Z\end{bmatrix} \succcurlyeq -\alpha \begin{bmatrix} I \\ Z\end{bmatrix}^\top N\begin{bmatrix} I \\ Z\end{bmatrix} = 0.
\]
Hence,~\ref{item:M2} holds. We complete the necessity part of the proof by noting that, by Theorem~\ref{thm:nsfl},~\ref{item:NS1} implies~\ref{item:NS3}. 

{\bf Sufficiency:} By Theorem~\ref{thm:nsfl},~\ref{item:NS1} holds if and only if~\ref{item:NS2} and~\ref{item:NS3} hold. Hence, it suffices here to show that~\ref{item:M2} and~\ref{item:NS3} imply~\ref{item:NS2}. Suppose that~\ref{item:M2} and~\ref{item:NS3} hold. By assumption, $Z\coloneqq -N_{22}^+N_{12}^\top$ satisfies
\begin{equation}
    \begin{bmatrix} I \\ Z\end{bmatrix}^\top N\begin{bmatrix}I \\ Z\end{bmatrix} = N_{11}-N_{12}N_{22}^+N_{12}^\top =0.\label{eq:IZNIZ}
\end{equation}
First, we show that $\ker N_{22}\subset\ker M_{22}$ and $\ker N_{22}\subset\ker M_{12}$. By $(I-N_{22}N_{22}^+)N_{12}^\top = 0$, we have $\im N_{12}^\top \subset \im N_{22}$ and, equivalently, $\ker N_{22}\subset \ker N_{12}$. Let $\xi\in\ker N_{22}$. By~\eqref{eq:IZNIZ}, it follows that 
\[
    \begin{bmatrix} I \\ Z + \xi\eta^\top\end{bmatrix}^\top N\begin{bmatrix} I \\ Z+ \xi\eta^\top\end{bmatrix}=0
\]
for any $\eta\in\mathbb{R}^{n}$. Thus, by~\ref{item:M2}, we have 
\[
    \begin{bmatrix} I \\ Z+\xi\eta^\top\end{bmatrix}^\top M\begin{bmatrix}I \\ Z+\xi\eta^\top\end{bmatrix}\succcurlyeq 0.
\]
Since $M_{22}\preccurlyeq 0$, it must hold that $\ker N_{22}\subset\ker M_{22}$. It follows that
\begin{equation}\label{eq:M12ker}
    \begin{bmatrix}
        I \\ Z
    \end{bmatrix}^\top M\begin{bmatrix}
         I \\ Z
        \end{bmatrix} + M_{12}\xi\eta^\top + \eta\xi^\top M_{12}^\top\succcurlyeq 0,
\end{equation}
for any $\eta\in\mathbb{R}^{m}$. Take $\xi\in\ker N_{22}$. Since we have~\eqref{eq:M12ker} for any $\eta \in \mathbb{R}^{n}$, it must hold that $M_{12}\xi\eta^\top+\eta\xi^\top M_{12}^\top=0$ for all $\eta\in\mathbb{R}^{n}$. Hence, $M_{12}\xi\eta^\top$ is skew-symmetric for any $\eta\in\mathbb{R}^{n}$. Let $\eta = (\eta_1,\eta_2,\hdots,\eta_n)$ and $v=(v_1,v_2,\hdots,v_n)=M_{12}\xi$. Then, for any $\eta\in\mathbb{R}^{n}$, the diagonal elements of $M_{12}\xi\eta^\top=v\xi^\top$ must satisfy $v_i\eta_i=0$ due to the skew-symmetry. It follows that $v_i=0$ for all $i\in\{1,2,\hdots,n\}$ and, thus, $v=M_{12}\xi=0$. We conclude that $\xi\in\ker M_{12}$ and, thus, $\ker N_{22}\subset\ker M_{12}$.

Next, we show that~\ref{item:NS2} holds. To this end, let $x\in\mathbb{R}^{n+m}$ be such that $x^\top Nx=0$, and let $T = \begin{bmatrix}\begin{smallmatrix}I & 0\\Z & I\end{smallmatrix}\end{bmatrix}$, such that $x=Tq$ for some $q=(q_1,q_2)\in\mathbb{R}^{n+m}$. Since $x^\top Nx=0$, we find, using $N_{11}-N_{12}N_{22}^+N_{12}^\top = 0$ and $(I-N_{22}N_{22}^+)N_{12}^\top = 0$, that $q^\top T^\top NTq = q_2^\top N_{22}q_2=0$.
Since $N_{22}\preccurlyeq 0$ by assumption, it must hold that $q_2\in\ker N_{22}$ and, thus, $q_2\in\ker M_{12}$ and $q_2\in\ker M_{22}$. It follows that
\begin{equation}\label{eq:q1M11q1}
    x^\top Mx = q_1^\top \begin{bmatrix} I \\ Z\end{bmatrix}^\top M\begin{bmatrix}I \\ Z\end{bmatrix} q_1\stackrel{\ref{item:M2}}{\geqslant}0.
\end{equation}
Thus,~\ref{item:M2} and~\ref{item:NS3} indeed imply~\ref{item:NS2}.\hfill\qed

\end{document}